\documentclass[a4paper,12pt]{article}
\usepackage{amssymb,amsmath,amsthm,latexsym}
\usepackage{amsfonts}
	\usepackage{amsfonts}
\usepackage{graphicx}
\usepackage[pdftex,bookmarks,colorlinks=false]{hyperref}
\usepackage{verbatim}
\usepackage{caption}
\usepackage{subcaption}
\usepackage[hmargin=1.15in,vmargin=1.15in]{geometry}
\usepackage{fancyhdr}
\newtheorem{theorem}{Theorem}[section]

\newtheorem{corollary}[theorem] {Corollary}
\newtheorem{definition}[theorem]{Definition}

\newtheorem{problem}{Problem}
\newtheorem{proposition}[theorem]{Proposition}

\title{\bf Arithmetic Intger Additive Set-Indexers of Graph Operations}
\author{N K Sudev \\ {\small Department of Mathematics} \\ {\small Vidya Academy of Science \& Technology}\\  {\small Thrissur - 680501, India.}\\ {\small email: {\em sudevnk@gmail.com}} \\ \\ K A Germina \\ {\small Department of Mathematics}\\ {\small School of Mathematical \& Physical Sciences}\\{\small Central University of Kerala} \\ {\small  Kasaragod-671316, India.}\\ {\small email:{\em srgerminaka@gmail.com}}}
\date{}
\begin{document}
\maketitle

\begin{abstract}
An integer additive set-indexer is an injective function $f:V(G)\to 2^{\mathbb{N}_0}$ such that the induced function $g_f:E(G) \to 2^{\mathbb{N}_0}$ defined by $g_f (uv) = f(u)+ f(v)$ is also injective. A graph $G$ which admits an IASI is called an IASI graph. An arithmetic integer additive set-indexer is an integer additive set-indexer $f$, under which the set-labels of all elements of a given graph $G$ are arithmetic progressions. In this paper, we discuss about admissibility of arithmetic integer additive set-indexers by certain graph operations and certain products of graphs. 
\end{abstract}
\textbf{Key words}: Integer additive set-indexers, arithmetic integer additive set-indexers, isoarithmetic integer additive set-indexers,biarithmetic integer additive set-indexers.\\
\textbf{AMS Subject Classification : 05C78} 

\section{Introduction}

For all  terms and definitions, not defined specifically in this paper, we refer to \cite{FH} and for more about graph labeling, we refer to \cite{JAG}. Unless mentioned otherwise, all graphs considered here are simple, finite and have no isolated vertices. All sets mentioned in this paper are finite sets of non-negative integers. We denote the cardinality of a set $A$ by $|A|$.

\begin{definition}\label{D2}{\rm
\cite{GA} An {\em integer additive set-indexer} (IASI, in short) is defined as an injective function $f:V(G)\rightarrow 2^{\mathbb{N}_0}$ such that the induced function $g_f:E(G) \rightarrow 2^{\mathbb{N}_0}$ defined by $g_f (uv) = f(u)+ f(v)$ is also injective.  A graph $G$ which admits an IASI is called an IASI graph.}
\end{definition}

\begin{definition}\label{D3}{\rm
The cardinality of the labeling set of an element (vertex or edge) of a graph $G$ is called the {\em set-indexing number} of that element.}
\end{definition}

In \cite{GS1}, the vertex set $V$ of a graph $G$ is defined to be {\em $l$-uniformly set-indexed}, if all the vertices of $G$ have the set-indexing number $l$.

By the term, an arithmetically progressive set, (AP-set, in short), we mean a set whose elements are in arithmetic progression. We call the common difference of the set-label of an element of the given graph, the {\em deterministic index} of that element.

\begin{proposition}\label{P-AIASI-d}
Let $f$ be a vertex-arithmetic IASI defined on $G$. If the set-labels of vertices of $G$ are AP-sets with the same common difference $d$, then $f$ is also an edge-arithmetic IASI of $G$.
\end{proposition}

\begin{definition}{\rm
\cite{GS7} An {\em arithmetic integer additive set-indexer} is an integer additive set-indexer $f$, under which the set-labels of all elements of a given graph $G$ are the sets whose elements are in arithmetic progressions. A graph that admits an arithmetic IASI is called an {\em arithmetic IASI graph}.

If all vertices of $G$ are labeled by the set consisting of arithmetic progressions, but the set-labels are not arithmetic progressions, then the corresponding IASI may be called {\em semi-arithmetic IASI}.}
\end{definition}

\begin{theorem}\label{T-AIASI-g}
\cite{GS7} A graph $G$ admits an arithmetic IASI graph $G$ if and only if for any two adjacent vertices in $G$, the deterministic index of one vertex is a positive integral multiple of the deterministic index of the other vertex and this positive integer is less than or equal to the cardinality of the set-label of the latter vertex.
\end{theorem}

\begin{definition}{\rm
\cite{GS8} If all the set-labels of all elements of a graph $G$ consist of arithmetic progressions with the same common difference $d$, then the corresponding IASI is called {\em isoarithmetic IASI}. }
\end{definition}

\begin{definition}{\rm
\cite{GS9} Let $f$ be an arithmetic IASI of a graph $G$. For two vertices $v_i$ and $v_j$ of $G$, let the common differences of $f(v_i)$ and $f(v_j)$ be $d_i$ and $d_j$ respectively.  If either of $d_i$ and $d_j$, for the adjacent vertices $v_i$ and $v_j$, is a positive integral multiple of the other, then $f$ is called a {\em biarithmetic IASI}. For $k\in \mathbb{N}_0$, if $d_i=k.d_j$ for all adjacent vertices $v_i$ and $v_j$ in $G$, then $f$ is called a {\em biarithmetic IASI} of $G$. 

If the value of $k$ is unique for all pairs of adjacent vertices of a biarithmetic IASI graph $G$, then that biarithmetic IASI is called {\em identical biarithmetic IASI} and $G$ is called an {\em identical biarithmetic IASI graph}.}
\end{definition}

\begin{theorem}\label{T-AIASI-SG}
\cite{GS7,GS8,GS9} A subgraph of an arithmetic (or isoarithmetic  or biarithmetic) IASI graph is also an arithmetic IASI graph (or isoarithmetic or biarithmetic) IASI graph.
\end{theorem}

\begin{theorem}\label{T-AIIBA1}
\cite{GS9} A graph $G$ admits an identical biarithmetic IASI if and only if it is bipartite.
\end{theorem}

\newpage

\section{New Results}

In this paper, we investigate the admissibility of arithmetic integer additive set-indexers by different operations and certain products of arithmetic IASI graphs.

\subsection{Arithmetic IASIs of Graph Operations}

The following result establishes the admissibility of the union of two arithmetic IASI graphs.

\begin{proposition}\label{T-IIASI-GJ1}
The union of two arithmetic IASI graphs admits an arithmetic IASI graph.
\end{proposition}
\begin{proof}
let $f_1$ and $f_2$ be the arithmetic IASIs defined on $G_1$ and $G_2$ respectively. Define a function $f$ on $G=G_1\cup G_2$ by
	\begin{equation*}
	f(v)=
	\begin{cases}
	f_1(v) & \text{if}~~v\in G_1\\
	f_2(v) & \text{if}~~v\in G_2.
	\end{cases}
	\end{equation*}
Since both $f_1$ and $f_2$ are arithmetic IASIs, then $f$ is also an arithmetic IASI on $G_1\cup G_2$.
\end{proof}

\begin{proposition}\label{T-IIASI-GJ1a}
The union of two isoarithmetic IASI graphs admits an isoarithmetic IASI graph if and only if all the vertices in both $G_1$ and $G_2$ have the same deterministic index.
\end{proposition}
\begin{proof}
Let $f_1$ and $f_2$ be the isoarithmetic IASIs defined on $G_1$ and $G_2$ respectively. Define a function $f$ on $G=G_1\cup G_2$ by
	\begin{equation*}
	f(v)=
	\begin{cases}
	f_1(v) & \text{if}~~v\in G_1\\
	f_2(v) & \text{if}~~v\in G_2.
	\end{cases}
	\end{equation*}
Assume that the vertices of both and $G_1$ and $G_2$ have the same deterministic index. Then, all the vertices of $G_1\cup G_2$ have the same deterministic index. Therefore, $f$ is  arithmetic IASI on $G_1\cup G_2$.

Conversely, assume that $G_1\cup G_2$ admits an isoarithmetic IASI, say $f$. Therefore, by Theorem \ref{T-AIASI-SG}, its subgraphs $G_1$ and $G_2$ also admit isoarithmetic IASIs which are the restrictions $f_1$ and $f_2$ of $f$ to $G_1$ and $G_2$ respectively.
\end{proof}

%\subsection{Arithmetic IASIs of Graph Joins}

\noindent Now, recall the definition of the join of two graphs.
\begin{definition}{\rm
\cite{FH} Let $G_1(V_1,E_1)$ and $G_2(V_2,E_2)$ be two graphs. Then, their {\em join} (or {\em sum}), denoted by $G_1+G_2$, is the graph whose vertex set is $V_1\cup V_2$ and edge set is $E_1\cup E_2\cup E_{ij}$, where $E_{ij}=\{u_iv_j:u_i\in G_1,v_j\in G_2\}$. }
\end{definition}

\begin{theorem}\label{T-IIASI-GJ2}
The join of two arithmetic IASI graphs admits an arithmetic IASI if and only if the deterministic index of every vertex of one graph is an integral multiple or divisor of the deterministic index of every vertex of the other graph, where this integer is less than or equal to the set-indexing number of the vertex having smaller deterministic index.
\end{theorem}
\begin{proof}
Let $G_1$ and $G_2$ be the given arithmetic IASI graphs. Let $E_{ij}=\{u_iv_j:u_i\in G_1,v_j\in G_2\}$ so that $G_1+G_2=G_1\cup G_2 \cup \langle E_{ij} \rangle$.

Assume that $G_1+G_2$ admits an arithmetic IASI, say $f$. Therefore, by Theorem \ref{T-AIASI-g}, for all edges in $\langle E_{ij} \rangle$ also, the deterministic index of one end vertex is an integral multiple of the deterministic index of the other end vertex, where this integer is less than or equal to the cardinality of the set-label of the latter vertex.  Since every vertex, say $u_i$, in $G_1$ is adjacent to every vertex, say $v_j$ of $G_2$ (and vice versa),  deterministic index of $u_i$ is either a multiple or a divisor of the deterministic index of $v_j$, where this integer is less than or equal to the set-indexing number of the vertex having smaller deterministic index.

Conversely, assume, without loss of generality, that the deterministic index of every vertex of $G_1$ is a multiple or a divisor of the deterministic index of every vertex of $G_2$ such that this integer is less than or equal to the set-indexing number of the vertex having smaller deterministic index. Hence, for every edge in $G_1+G_2$, the deterministic index of one end vertex is an integral multiple of the deterministic index of the other end vertex, since both $G_1$ and $G_2$ are arithmetic IASI graph. Then, by Theorem \ref{T-AIASI-g}, $G_1+G_2$ admits an arithmetic IASI.
\end{proof}

The following results establish the admissibility of arithmetic and isoarithmetic IASIs by the join of two isoarithmetic IASI graphs.

\begin{proposition}
The join of two isoarithmetic IASI graphs is an arithmetic IASI graph if and only if the deterministic index of the elements of one graph is an integral multiple of the deterministic index of the elements of the other, where this integer is less than or equal to the minimum among cardinalities of all set-labels of the former graph.
\end{proposition}
\begin{proof}
Let $G_1$ and $G_2$ admit isoarithmetic IASIs $f_1$ and $f_2$ respectively and let $E_{ij}=\{u_iu_j:u_i\in G_1,u_j\in G_2\}$ be such that $G_1+G_2=G_1\cup G_2 \cup \langle E_{ij} \rangle$. Note that all the elements of $G_1$ have the same deterministic index, say $d_1$ and all the elements of $G_2$ also have the same deterministic index, say $d_2$. Let that $d_1\ne d_2$.

Now, assume that the sum of two isoarithmetic IASI graphs is an arithmetic IASI graph. Then, by Theorem \ref{T-AIASI-g}, for all edges in $E_{ij}$, the deterministic index of one end vertex is an integral multiple of the deterministic index of the other, with this integer is less than or equal to the cardinality of the set-label of the former vertex.  Without loss of generality, let $d'_j=k~d_i$, where $d_i$ is the deterministic index of the vertex $v_i$ in $G_1$ and $d'_j$ is the deterministic index of the vertex $u_j$ in $G_2$ and $k$ is a positive integer $|f_1(v_i)|$. Therefore, the deterministic index of the elements of $G_2$ is an integral multiple of the deterministic index of the elements of $G_1$, where this integer is less than or equal to the minimum among the cardinalities of all set-labels of the $G_1$.

Conversely, assume, without loss of generality, that the deterministic index of the elements of $G_2$ is an integral multiple of the deterministic index of the elements of $G_1$, where this integer is less than or equal to the minimum among the cardinalities of all set-labels of $G_1$.  Therefore, by Theorem \ref{T-AIASI-g}, $G_1+G_2$ admits an arithmetic IASI.
\end{proof}

\begin{proposition}\label{P-IIASI-GJ3}
The join of two isoarithmetic IASI graphs admits an isoarithmetic IASI if and only if all the vertices in both $G_1$ and $G_2$ have the same deterministic index.
\end{proposition}
\begin{proof}
First assume that $G_1+G_2$ admits an isoarithmetic IASI, say $f$. Then, the deterministic index of all vertices of $G_1+G_2$ is the same,say $d$. Since, $G_1$ and $G_2$ are subgraphs of $G_1+G_2$, by Theorem \ref{T-AIASI-SG}, $G_1$ and $G_2$ admits the induced isoarithmetic IASI of $f$. All the vertices in both $G_1$ and $G_2$ have the same deterministic index. 

Conversely, the vertices in both $G_1$ and $G_2$ have the same deterministic index. Let $E_{ij}=\{u_iv_j:u_i\in G_1,v_j\in G_2\}$. Then, for any edge in $\langle E_{ij} \rangle$ must the deterministic index its end vertices are the same as that of both $G_1$ and $G_2$. Hence, by Theorem \ref{T-AIASI-g}, $G_1+G_2$ admits an isoarithmetic IASI.
\end{proof}

\begin{proposition}
The join of two isoarithmetic IASI graphs is an arithmetic IASI graph if and only if the deterministic index of the elements of one graph is an integral multiple of the deterministic index of the elements of the other.
\end{proposition}

\begin{proposition}\label{C-BIASI-GJ4}
Let $G_1$ and $G_2$ be two graphs which admit identical biarithmetic IASI. Then, $G_1+G_2$ does not admit an identical biarithmetic IASI.
\end{proposition}
\begin{proof}
If possible, let $G_1+G_2$ admits an identical biarithmetic IASI. Since $G_1$ and $G_2$ are identical biarithmetic IASI graphs, for every pair of adjacent vertices in them, the deterministic index of one is a positive integral multiple of the deterministic index of the other and this positive integer is unique for all such pair vertices in $G_1$ and $G_2$. Let $v_i$ be a vertex of $G_1$ with deterministic index $ d_i$ and let $u_j$ and $u_l$ be two adjacent vertices in $G_2$ with deterministic indices $d'_j$ and $d'_l$ respectively. By Theorem \ref{T-AIASI-g}, $d'_l=k.d'_j$ for some positive integer $k$.   

 Now, $v_i$ is adjacent to both and $u_j$ and $u_l$ in $G_1+G_2$. Then, by Theorem \ref{T-AIASI-g}, we have $d_i=k.d'_j$, $d'_l=k.d'_j$ and $d_i=k.d'_j$, all of which cannot hold simultaneously. Therefore, $G_1+G_2$ does not admit an identical biarithmetic IASI. Hence, $G_1+G_2$ does not admit an identical biarithmetic IASI.
\end{proof}

\subsection{Arithmetic IASIs of Graph Products}

We discuss the admissibility of arithmetic IASI by certain graph products. First, recall the definition of the cartesian product of two graphs. 

\begin{definition}{\rm 
\cite{FH} Let $G_1(V_1,E_1)$ and $G_2(V_2,E_2)$ be two graphs.Then, the {\em cartesian product} of $G_1$ and $G_2$, denoted by $G_1\times G_2$, is the graph with vertex set $V_1\times V_2$  defined as follows. Let $u=(u_1, u_2)$ and $v=(v_1,v_2)$ be two points in $V_1\times V_2$. Then, $u$ and $v$ are adjacent in $G_1\times G_2$ whenever [$u_1=v_1$ and $u_2$ is adjacent to $v_2$] or [$u_2=v_2$ and $u_1$ is adjacent to $v_1$]. If $|V_i|=p_i$ and $E_i=q_i$ for $i=1,2$., then $|V(G_1\times G_2)|=p_1p_2$  and $i=1,2$ and $|E(G_1\times G_2)|=p_1q_2+p_2q_1$.}
\end{definition}

The cartesian product $G_1\times G_2$ may be viewed as follows. Make $p_2$ copies of $G_1$. Denote these copies by $G_{1_i}$, which corresponds to the vertex $v_i$ of $G_2$. Now, join the corresponding vertices of two copies $G_{1i}$ and $G_{1j}$ if the corresponding vertices $v_i$ and $v_j$ are adjacent in $G_2$. Thus, we view the product $G_1\times G_2$ as a union of $p_2$ copies of $G_1$ and a finite number of edges connecting two copies $G_{1i}$ and $G_{1j}$ of $G_1$ according to the adjacency of the corresponding vertices $v_i$ and $v_j$ in $G_2$, where $1\le i\neq j\le p_2$.

The following theorem establishes the admissibility of arithmetic IASI by the cartesian product of two arithmetic IASI graphs. 

\begin{theorem}\label{T-IIASI-GJ5}
The cartesian product of two arithmetic IASI graphs $G_1$ and $G_2$ admits an arithmetic IASI if and only if, for any pair of corresponding vertices of the copies of $G_1$ (or $G_2$) which are adjacent in $G_1\times G_2$, the deterministic index of one vertex is an integral multiple (or a divisor) of the deterministic index of the  other vertex, where this integer is less than or equal to the set-indexing number of the vertex having smaller deterministic index.
\end{theorem}
\begin{proof}
Let $U=\{u_1,u_2,u_3,\ldots, u_n\}$ be the vertex set of $G_1$ and $V=\{v_1,v_2,v_3,\ldots, v_m\}$ be the vertex set of $G_2$. Let $G_{1j};~ 1\le j\le m$, be the $m$ copies of $G_1$ in $G_1\times G_2$. Therefore, $G_{1j}= \langle U_j \rangle$ where $U_j=\{u_{ij}:1\le i \le |E(G_1)|, 1\le j \le |E(G_2)|\}$. Now, for all values of $j$, the graphs induced by the set of vertices $\{u_{ij}:1\le i \le |E(G_1)|\}$ are graphs isomorphic of $G_1$ and similarly for all values of $i$, the graphs induced by the set of vertices $\{u_{ij}:1\le j \le |E(G_2)|\}$ are the graphs isomorphic of $G_2$. Without loss of generality, let $\langle U_1 \rangle = G_1$ and $\langle V_1 \rangle = G_2$, where $V_1=\{u_{i1}:1\le i\le |E(G_2)|\}$.  Also, the corresponding vertices of $G_{1r}$ and $G_{1s}$ are adjacent in $G_1\times G_2$ if the vertices $v_r$ and $v_s$ adjacent in $G_2$. 

Now, let $f_1$ and $f_2$ be the arithmetic IASIs of $G_1$ and $G_2$ respectively. Since $G_1$ is an arithmetic IASI graph, for two adjacent vertices $u_r$ and $u_s$ in $G_1$, we have $d_s=k_l.d_r$, where $k_l \le |f_1(u_r)|$, is a positive integer with $1\le l \le |E(G_1)|$. Similarly, Since $G_2$ is an arithmetic IASI graph, for two adjacent vertices $v_r$ and $v_s$ in $G_2$, we have $d'_s=k'_l.d'_r$, where $k'_l \le |f_2(v_r)|$, is a positive integer with $1\le l \le |E(G_2)|$, where $d_i$ is the deterministic index of the vertex $u_i$ in $G_1$ and $d'_j$ is the deterministic index of the vertex $v_j$ in $G_2$. Label the vertices of $\langle U_1 \rangle$ by the same set-labels of $G_1$ itself and label the vertices of the copies $\langle U_r \rangle$ and $\langle U_s \rangle$ in such a way that the deterministic indices of all vertices in $U_s$ are integral multiple of the deterministic indices of the corresponding vertices of $U_r$ by a unique positive integer $k'_l$, $k'_l$ being the same positive integer given by $k'_l=\frac{d'_s}{d'_r}$ where $d'_r$ and $d'_s$ are the deterministic indices of the vertices $v_r$ and $v_s$ in $G_2$ corresponding to the copies $U_r$ and $U_s$ respectively. Then, for every pair of adjacent vertices in $G_1\times G_2$, the deterministic index of one vertex is an integral multiple of the deterministic index of the other. Hence, by Theorem \ref{T-AIASI-g}, $G_1\times G_2$ is an arithmetic IASI graph.

Conversely, assume that $G_1\times G_2$ admits an arithmetic IASI. Now we can take the same set-labels of the copy $U_1$ as the set-labels of the vertices of $G_1$ and the same set-labels of the graph $\langle V_1\rangle$, defined above, as the set-labels of the graph $G_2$. Since the set-labels of all elements of $G_1\times G_2$ are AP-sets, these set-labelings of $G_1$ and $G_2$ are arithmetic IASIs. Therefore $G_1$ and $G_2$ are arithmetic IASI graphs.
\end{proof}

Invoking Proposition \ref{T-IIASI-GJ5}, we now establish the following results.

\begin{corollary}\label{C-IIASI-GJ7}
The cartesian product of two isoarithmetic IASI graphs admits an isoarithmetic IASI if and only if all vertices in both $G_1$ and $G_2$ hve the same deterministic index.
\end{corollary}
\begin{proof}
The proof is immediate form Theorem \ref{T-IIASI-GJ5} by taking $k_l=1$ and $k'_l=1$.
\end{proof}

Now, we investigate the admissibility of an arithmetic IASI by the cartesian product of two identical biarithmetic IASI graphs. %First, recall the following theorem on the identical biarithmetic IASI graphs, proved in \cite{GS9}.

Invoking Theorem \ref{T-AIIBA1} and Theorem \ref{T-IIASI-GJ5} we establish the following result.

\begin{proposition}\label{P-IIASI-GJ6}
The cartesian product of two identical biarithmetic IASI graphs admits an identical biarithmetic IASI. 
\end{proposition}
\begin{proof}
Since $G_1$ and $G_2$ are identical biarithmetic IASI graphs, by Theorem \ref{T-AIIBA1}, both $G_1$ and $G_2$ are bipartite. Therefore, since the cartesian product of two bipartite graphs is also a bipartite graph, $G_1\times G_2$ is bipartite. Hence, $G_1\times G_2$ admits an identical biarithmetic IASI. The labeling of the vertices in $G_1\times G_2$ follows from Theorem \ref{T-IIASI-GJ5} by taking $k_l=k'_l=k$, a unique positive integer.
\end{proof}

Next, we proceed to verify the admissibility of arithmetic IASI by the corona of two graphs. Now, recall the definition of corona of two graphs.

\begin{definition}{\rm
\cite{FH} By the term {\em corona} of two graphs $G_1$ and $G_2$, denoted by $G_1\circ G_2$, is the graph obtained taking one copy of $G_1$ (which has $p_1$ vertices) and $p_1$ copies of $G_2$ and then joining the $i$-th point of $G_1$ to every point in the $i$-th copy of $G_2$.}
\end{definition}

The number of vertices and edges in $G_1\circ G_2$ are $p_1(1+p_2)$ and $q_1+p_1q_2+p_1p_2$ respectively, where $p_i$ and $q_i$ are the number of vertices and edges of the graph $G_i, i=1,2$.

\begin{theorem}\label{T-IIASI-GJ8}
Let $G_1$ and $G_2$ are an arithmetic IASI graphs. Then, the corona $G_1\circ G_2$ admits an arithmetic IASI if and only if the deterministic index of every vertex of one graph is an integral multiple or a divisor of the deterministic index of every vertex of the other, where this integer is less than or equal to the set-indexing number of the vertex having smaler deterministic index..
\end{theorem}
\begin{proof}
Let $U=\{u_1,u_2,u_3,\ldots, u_n\}$ be the vertex set of $G_1$ and $V=\{v_1,v_2,v_3,\ldots, v_m\}$ be the vertex set of $G_2$. Let $G_{2i};~ 1\le i\le n$, be the $n$ copies of $G_2$ in $G_1\circ G_2$. Therefore, $G_{2i}= \langle V_i \rangle$ where $V_i=\{v_{ij}:1\le i \le |E(G_1)|~\text{and}~ 1\le j \le |E(G_2)|\}$.  Now, for any value of $i=r$, the graphs induced by the set of vertices $\{v_{rj}:1\le j \le |E(G_2)|\}$ is a graph isomorphic of $G_2$.  Without loss of generality, let $\langle V_1 \rangle = G_2$.   Also, all the vertices of $\langle V_i \rangle $ are adjacent to the vertex $u_i$ of $G_1$ in $G_1\circ G_2$.  

Now, assume that $G_1$ and $G_2$ admit arithmetic IASIs $f_1$ and $f_2$ respectively. Since $G_1$ is an arithmetic IASI graph, for two adjacent vertices $u_r$ and $u_s$ in $G_1$, we have $d_s=k_l.d_r$, where $k_l\le |f_1(u_r)|$ is a positive integer with $1\le l \le |E(G_1)|$. Similarly, since $G_2$ is an arithmetic IASI graph, for two adjacent vertices $v_r$ and $v_s$ in $G_2$, we have $d'_s=k'_l.d'_r$, where $k'_l\le |f_2(v_r)|$ is a positive integer with $1\le l \le |E(G_2)|$, where $d_i$ is the deterministic index of the vertex $u_i$ in $G_1$ and $d'_j$ is the deterministic index of the vertex $v_j$ in $G_2$. Label the vertices of $\langle V_1 \rangle$ by the same set-labels of $G_2$ itself and label the vertices of the copies $\langle V_r \rangle; ~ 1<r\le n$, by distinct sets in such a way that the deterministic indices of the corresponding vertices in $V_1$ and $V_r$ have the same deterministic indices.  Then, for every pair of adjacent vertices in $G_1\circ G_2$, the deterministic index of one vertex is an integral multiple of the deterministic index of the other. Hence, by Theorem \ref{T-AIASI-g}, $G_1\circ G_2$ is an arithmetic IASI graph.

Conversely, assume that $G_1\circ G_2$ admits an arithmetic IASI. Since $G_1$ is a subgraph of $G_1\circ G_2$, by Theorem \ref{T-AIASI-SG}, $G_1$ admits an arithmetic IASI. Also, we can take the same set-labels of the copy $V_1$ as the set-labels of the vertices of $G_2$ and the same set-labels of the component graph $\langle V_1\rangle$ of $G_1\circ G_2$, defined above, as the set-labels of the graph $G_2$. Since the set-labels of all elements of $G_1\times G_2$ are AP-sets, these set-labelings of $G_1$ and $G_2$ are arithmetic IASIs. Therefore $G_1$ and $G_2$ are arithmetic IASI graphs.
\end{proof}

Invoking Proposition \ref{T-IIASI-GJ5}, we now establish the following results.

\begin{proposition}\label{C-IIASI-GJ9}
The corona of two isoarithmetic IASI graphs admits an isoarithmetic IASI if and only if all vertices in both $G_1$ and $G_2$ hve the same deterministic index.
\end{proposition}
\begin{proof}
The proof is immediate form Theorem \ref{T-IIASI-GJ8}.
\end{proof}

\begin{proposition}\label{C-IIASI-GJ10}
The corona of two identical biarithmetic IASI graphs does not admit an identical biarithmetic IASI.
\end{proposition}
\begin{proof}
Since the corona of two bipartite graphs is not a bipartite graph, by theorem \ref{T-AIIBA1}, $G_1\circ G_2$ does not admit an identical biarithmetic IASI.
\end{proof}

\subsection{Arithmetic IASIs of Graph Complements}

In this section, we discuss the admissibility  of arithmetic  IASI by complements of given arithmetic IASI graphs. Note that the vertices of a graph $G$ and its complements have the same set-labels and hence the same deterministic indices. 

\begin{theorem}
The complement of an arithmetic IASI graph $G$ admits an arithmetic IASI if and only if the deterministic index of any vertex of $G$ is an integral multiple or divisor of the deterministic index of every other vertex of $G$.
\end{theorem}
\begin{proof}
First assume that the deterministic index of any vertex of $G$ is an integral multiple or divisor of the deterministic index of every other vertex of $G$. Hence, for any pair of adjacent vertices in $\bar{G}$ also, the deterministic index of one vertex is an integral multiple of the deterministic index of the other. Hence by Theorem \ref{T-AIASI-g}, $\bar{G}$ admits an arithmetic IASI.

Conversely, assume that $\bar{G}$ admits an arithmetic IASI. Since every pair of vertices in $V(G)$ are either adjacent in $G$ or in its complement $\bar{G}$ and both $G$ and $\bar{G}$ are arithmetic IASI graphs, for every pair of vertices the deterministic index of one vertex must be a multiple of the deterministic index of the other. 
\end{proof}

\begin{proposition}
The complement of an isoarithmetic IASI graph is also an isoarithmetic IASI graph.
\end{proposition}
\begin{proof}
Let $G$ be an isoarithmetic IASI graph. Then, the deterministic index of all vertices of $G$ (and $\bar{G}$) are the same. Therefore, $\bar{G}$ admits an isoarithmetic IASI.
\end{proof}

\begin{proposition}
The complement of an identical biarithmetic IASI never admits an identical biarithmetic IASI.
\end{proposition}
\begin{proof}
The complement of a bipartite graph $G$ is not a bipartite graph. Hence, by Theorem \ref{T-AIIBA1}, $\bar{G}$ does not admit an identical biarithmetic IASI.
\end{proof}

\section{Conclusion}

In this paper, we have discussed some characteristics of certain graphs operations and products which admit arithmetic IASIs. We have not addressed certain Problems in this area which are still open.

The following are some of the open problems we have identified in this area.

\begin{problem}{\rm 
Examine the necessary and sufficient conditions for the existence of non-identical biarithmetic IASIs for given graphs.}
\end{problem} 

\begin{problem}{\rm 
Discuss the necessary and sufficient conditions for the existence of arithmetic IASIs for some other products, such as lexicographic product, tensor product, strong product, rooted product etc., of arithmetic IASI graphs.}
\end{problem} 

\begin{problem}{\rm 
Discuss the necessary and sufficient conditions for the existence of arithmetic IASIs of all types for certain powers of arithmetic IASI graphs.}
\end{problem}

\begin{problem}{\rm 
Discuss the necessary and sufficient conditions for the existence of arithmetic IASIs for different graph classes having arithmetic IASIs.}
\end{problem}

\begin{problem}{\rm 
Characterise certain graphs and graph classes in accordance with their admissibility of identical and non-identical biarithetic IASIs.}
\end{problem}

The IASIs under which the vertices of a given graph are labeled by different standard sequences of non negative integers, are also worth studying.   The problems of establishing the necessary and sufficient conditions for various graphs and graph classes to have certain IASIs still remain unsettled.


\begin{thebibliography}{15}
\bibitem {A10} B D Acharya, (1990). {\em Arithmetic Graphs}, J. Graph Theory, {\bf 14}(3), 275-299. 
\bibitem {AGA} B D Acharya, K A Germina and T M K Anandavally, {\em Some New Perspective on Arithmetic Graphs} In {\bf Labeling of Discrete Structures and Applications}, (Eds.: B D Acharya, S Arumugam and A Rosa), Narosa Publishing House, New Delhi, (2008), 41-46.
\bibitem {BM1} J A Bondy and U S R Murty, (1976). {\bf Graph Theory with Apllications}, North-Holland, Amsterdam.
\bibitem {CZ} G Chartrand and P Zhang, (2005). {\bf Introduction to Graph Theory}, McGraw-Hill Inc.
\bibitem{ND} N Deo, (1974). {\bf Graph Theory with Applications to Engineering and Computer Science}, Prentice-Hall.
\bibitem {FFH} R Frucht and F Harary (1970). {\em On the Corona of Two Graphs}, Aequationes Math., {\bf 4}(3), 322-325.
\bibitem {JAG} J A Gallian, (2011). {\em A Dynamic Survey of Graph Labelling}, The Electronic Journal of Combinatorics (DS 16).
\bibitem {GA} K A Germina and T M K Anandavally, (2012). {\em Integer Additive Set-Indexers of a Graph:Sum Square Graphs}, Journal of Combinatorics, Information and System Sciences, {\bf 37}(2-4), 345-358.
\bibitem {GS1} K A Germina, N K Sudev, (2013). {\em On Weakly Uniform Integer Additive Set-Indexers of Graphs}, Int. Math. Forum., {\bf 8}(37), 1827-1834.
%\bibitem {GS2} K A Germina, N K Sudev, {\em Some New Results on Strong Integer Additive Set-Indexers}, Communicated.
\bibitem {HIS} R Hammack, W Imrich and S Klavzar (2011). {\em Handbook of Product graphs}, CRC Press.
\bibitem {FH}  F Harary, (1969). {\bf Graph Theory}, Addison-Wesley Publishing Company Inc.
\bibitem {SMH} S M Hegde, (1989). {\em Numbered Graphs and Their Applications}, PhD Thesis, Delhi University.
\bibitem {IK} W Imrich, S Klavzar, (2000). {\bf Product Graphs: Structure and Recognition}, Wiley.
\bibitem {IKR} W Imrich, S Klavzar and D F Rall, (2008). {\bf Topics in Graph Theory: Graphs and Their Cartesian Products}, A K Peters.
\bibitem {KDJ} K D Joshi, {\bf Applied Discrete Structures}, New Age International, (2003).
\bibitem {MBN} M B Nathanson (1996). {\bf Additive Number Theory, Inverse Problems and Geometry of Sumsets}, Springer, New York.
%\bibitem {MBN2} M B Nathanson (2000). {\bf Elementary Methods in Number Theory}, Springer-Verlag, New York.
%\bibitem {MBN3} M B Nathanson (1996). {\bf Additive Number Theory: The Classical Bases}, Springer-Verlag, New York.
\bibitem {GS0} N K Sudev and K A Germina, (2014). {\em A Note on Integer Additive Set-Indexers of Graphs}, Int. J. Math. Sci.\& Engg. Applications, {\bf 8}(2), 11-22.
\bibitem {GS7} N K Sudev and K A Germina, {\em On Arithmetic Integer Additive Set-Indexers of Graphs}, submitted.
\bibitem {GS8} N K Sudev and K A Germina, {\em On Isoarithmetic Integer Additive Set-Indexers of Graphs}, submitted.
\bibitem {GS9} N K Sudev and K A Germina, {\em Biarithmetic Integer Additive Set-Indexers of Graphs}, submitted.
\bibitem {RJT} R J Trudeau, (1993). {\bf Introduction to Graph Theory}, Dover Pub., New York.
\bibitem {DBW} D B West, (2001). {\bf Introduction to Graph Theory}, Pearson Education Inc.

\end{thebibliography}
\end{document}